\documentclass[12pt]{article}
\usepackage{amsfonts, amsmath, amsthm}

\newtheorem{theorem}{Theorem}[section]
\newtheorem{lemma}[theorem]{Lemma}

\newtheorem{cor}[theorem]{Corollary}
\newtheorem{prop}[theorem]{Proposition}

\def\QQ{\mathbb{Q}}
\def\RR{\mathbb{R}}

\def\calO{\mathcal{O}}

\def\alg{\mathrm{alg}}
\def\an{\mathrm{an}}
\def\con{\mathrm{con}}

\def\imm{\mathrm{imm}}

\def\beq{\begin{equation}}
\def\eeq{\end{equation}}

\def\GK{\Gamma^K}
\def\GL{\Gamma^L}

\def\Gimm{\Gamma^{\imm}}
\def\Galg{\Gamma^{\alg}}

\def\Galgancon{\Galg_{\an, \con}}
\def\Galgcon{\Galg_{\con}}

\def\Gancon{\Gamma_{\an,\con}}
\def\Gcon{\Gamma_{\con}}
\def\Gimmancon{\Gimm_{\an,\con}}
\def\Gimmcon{\Gimm_{\con}}
\def\GKcon{\GK_{\con}}
\def\GLancon{\GL_{\an,\con}}
\def\GLcon{\GL_{\con}}

\DeclareMathOperator{\Gal}{Gal}

\begin{document}

\title{Quasi-unipotence of overconvergent $F$-crystals}
\author{Kiran S. Kedlaya}
\date{June 22, 2001}

\maketitle

\begin{abstract}
We prove that every overconvergent $F$-isocrystal over $k((t))$ is
quasi-unipotent (in the sense of Crew),
for $k$ a field of characteristic $p>0$.
\end{abstract}

Building on our earlier papers \cite{bib:me3}, \cite{bib:me4},
\cite{bib:me5},
we prove Crew's conjecture
that every overconvergent $F$-isocrystal over $k((t))$ is
quasi-unipotent.
Consequences of this result will be given in a subsequent paper.

Note that the result of this paper is what might be called the ``weak local
monodromy conjecture''; the ``strong'' version states that
any $(F, \nabla)$-module over the Robba ring is quasi-unipotent.
In this case, there is no generic
Newton polygon that can be used to construct the special Newton polygon.
It may be possible to refine the result of \cite{bib:me5} to allow this
case to be treated, but we do not presently see how to do so.

It is also worth pointing out that while quasi-unipotence depends on the 
presence both of the Frobenius and connection structure, our proof makes
minimal use of properties of $p$-adic differential equations.
Rather, we make a careful study of Frobenius structures over various rings.
These rings are carefully chosen, on the one hand
to allow various auxiliary constructions to be made, and on the other hand
to allow some of these constructions to be ``descended'' to the rings 
of actual interest.

By contrast, Yves Andr\'e (preprint) and Zoghman Mebkhout
(to appear) have announced proofs of
the strong local monodromy conjecture in case the residue field is the
algebraic closure of a finite field. Their approaches use results of Christol
and Mebkhout on the properties of $p$-adic differential equations and
minimal use of the Frobenius structure. The difference can be summarized
by saying that their proofs are ``$p$-adic analytic'' while ours is
``$t$-adic analytic''.

We retain the definitions and notations of our earlier papers
\cite{bib:me3}, \cite{bib:me4}, \cite{bib:me5}. In this notation,
the main result of this paper can be formulated as follows.
\begin{theorem}
  Let $M$ be an $(F, \nabla)$-crystal over $\GKcon$, with $\sigma =
\sigma_t$ standard. Then $M$ becomes
unipotent (as an $(F, \nabla)$-crystal)
over $\Gamma^L_{\an,\con} \otimes_{\calO} \calO'$
for some finite separable extension $L$
of $K$ and some finite extension $\calO'$ of $\calO$.
\end{theorem}
Recall that by an argument of Tsuzuki \cite[Proposition~3.1]{bib:me3}, this theorem implies
the corresponding statement for $\sigma$ arbitrary.

We add one notation: for $r>0$, $\Gamma^*_r$ (resp.\ $\Gamma^*_{\an,r}$)
 is the subring of
$\Gamma^*_{\con}$ (resp.\ $\Gamma^*_{\an,\con}$
consisting of $x = \sum x_i t^i$ with
$\limsup_{i \to -\infty} v_p(x_i)/(-i) \geq r$.

\section{Descending a slope filtration}

Recall the main result of \cite{bib:me4}.
\begin{theorem} \label{thm:constant}
Every $F$-crystal over $R = \Galgcon, \Gimmcon$
becomes constant over $R_{\an} \otimes_{\calO} \calO'$
for some finite extension $\calO'$ of $\calO$.
\end{theorem}
An alternate formulate of the result, avoiding the extension of scalars,
is that given an invertible matrix $A$ over $R$, there exist
$U$ over $R_{\an}$ such that $U^{-1}AU^\sigma$ has entries in $\calO$.

While this provides a new slope filtration over $\Galgcon$ of an
$F$-crystal over $\Gcon$, this slope filtration is not immediately useful
in studying $(F, \nabla)$-crystals, because $\Galgcon$ does not
carry a derivation. To use the filtration to prove quasi-unipotence,
we must first descend it to $\Gancon$ itself. This descent is the object
of this section.

Recall the norms $w_r$ defined on appropriate subrings of
$\Gamma^*_{\an,\con}$, and that the subring $\Omega^*_{\an,\con}$
is complete in the Fr\'echet topology given by these norms.
We extend these norms to matrices by declaring the norm of a matrix to be
the maximum norm of any of its entries.
Also, let $I_n$ denote the $n \times n$ identity matrix.

We say $x \in \Gimmancon$ is an \emph{$r$-semiunit} if 
there exists $c \in \calO$ and $i \in \QQ$ such that $w_r(x-ct^i) < w_r(x)$.
For $L$ a finite extension of $K$ (not necessarily separable), an
\emph{$r$-semiunit basis} is a sequence of $r$-semiunits $x_1, \dots, x_m$
whose reductions modulo $\pi$ form a basis for $L$ as a $K$-vector space.
\begin{lemma}
For any finite extension $L$ of $K$ (not necessarily separable),
there exists an extension $L_1$ of $L$ such that $\Gamma^{L_1}_{\con}$
admits an $r$-semiunit basis.
\end{lemma}
\begin{proof}
Let $x$ be an element of $\GLcon$ whose reduction modulo $\pi$ generates
$L$ over $K$. Then for $m$ sufficiently large, $x^{\sigma^{-m}}$ is an
$r$-semiunit. Thus if we take $L_1 = L^{1/p^m}$, we may use
powers of $x^{\sigma^{-m}}$ as an $r$-semiunit basis.
\end{proof}

Since $\Galgancon$ has the property that any finitely generated ideal
is principal \cite[Lemma~3.15]{bib:me5},
the following lemma can be deduced from
\cite[Theorem~III.7.9]{bib:lang}.
Recall that an \emph{elementary matrix} is a matrix obtained from the
identity by adding a multiple of one row to another, swapping two rows, or
multiplying one row by a unit.
\begin{lemma} \label{lem:elem}
For $r>0$,
any invertible matrix over $\Galg_{\an,r}$ can be factored as a product of
elementary matrices.
\end{lemma}
\begin{cor} \label{cor:approx}
Let $U$ be an invertible matrix over $\Galg_{\an,r}$ for some $r>0$.
Then there exists an
invertible matrix $V$ over $\GL_{\an,\con}$ 
such that $w_r(VU-I_n) < 1$
and $w_r((VU)^{-1} - I_n) < 1$.
\end{cor}
\begin{proof}
Factor $U^{-1}$ as a product of elementary matrices in $\Galg_{\an,r}$.
It suffices to note that for each elementary matrix $W$ in this product
and each $\epsilon > 0$,
there exists $X$ over $\GL_{\an,\con}$ for some finite extension $L$ 
of $K$ such that $w_r(X - U) < \epsilon$ and $w_r(X^{-1} - U^{-1}) < \epsilon$.
This is vacuous for permutation matrices and easily verified for the
other two kinds of matrices.
\end{proof}

\begin{lemma} \label{lem:anfact1}
Let $A$ be an $n \times n$ matrix over $\Gancon$, and
$D$ an invertible $n \times n$ diagonal matrix.
Suppose there exists an $n \times m$ matrix $U$ over $\Gimm_{\an,\con}$
such that $AU^{\sigma} = D U$.
Then there exists an invertible matrix $S$ over $\GLancon$,
for some finite extension $L$ of $K$,
such that $D^{-1}SAS^{-\sigma}$ has entries in $\GLcon$
and is congruent to $I_n$ modulo $\pi$.
\end{lemma}
\begin{proof}
Fix $r>0$ such that $w_r(U)$ and $w_r(U^\sigma)$
are defined, and let $v$ be the smallest positive valuation of
$\calO$ and $\Delta = \max_{i,j} v_p(D_{ii}/D_{jj})$.
By Corollary~\ref{cor:approx}, for any $\epsilon > 0$,
we can change basis over $\GL_{\an,r}$ for some
finite extension $L$ of $K$ to ensure that $w_r(U^\sigma-I_n) < \epsilon$
and $w_r(D^{-1} U D - I_n) < \epsilon$. In particular, choose $s$
large enough that $(p-1)rs - \max_{i,j} v_p(D_{ii}/D_{jj}) >v$ and
take $\epsilon = p^{-rs}$. We may also assume that $L$ contains
$vr$, and that $L$ has a $pr$-semiunit basis, by enlarging $L$ suitably.

We construct a convergent sequence $\{S_l\}_{l=0}^\infty$ over 
$\GLancon$ such that $w_r(D^{-1}S_l A S_l^{-\sigma} - I_n) < \epsilon$ for
each $l$, as follows. First set $S_0 = 1$. Given $S_l$,
if $D^{-1} S_l A S_l^{-\sigma}-I_n$ has no coefficients of norm at least
1,
the proof is complete and we stop. Otherwise,
write $D^{-1} S_l A S_l^{-\sigma} = \sum_i T_i t^i$.
If the smallest positive index $i$ for which $|T_i| > 1$
gives $|T_i| = 1$, let $j=j_l$ be the index, let $V$ be a matrix over
$\GL_{r}$ whose entries are $rp$-semiunits such that
$w_r(I_n + T_j t^j - V) < w_r(T_j)$, and put $W = V^{-1}$.
In this case we say $l$ is of the ``first type''.
Otherwise, let $j=j_l$ be the smallest positive index at which $r_l = 
\max (\log_p |T_j|)/j$ is realized.
By construction, $w_{r_l}(D^{-1}S_lAS_l^{-\sigma})
= w_{r_l}(S_l^\sigma A^{-1} S_l^{-1}D) = 1$. Thus the lowest
order terms of $S_l^{\sigma} A^{-1} S_l^{-1} D$ form an invertible matrix
over the ring of elements of $\GL_{\an,r}$ of $w_{r_l}$-norm
1 modulo those of norm less than 1. That ring is isomorphic to $k[u, u^{-1}]$,
and so is a principal ideal domain. Thus by Lemma~\ref{lem:elem}, the matrix
can be factored into elementary matrices. By lifting each of these into
an invertible matrix over $\GL_{\an,rp}$
(possible because $L$ admits an $rp$-semiunit basis)
and multiplying together, we get a matrix which we again call $W$.
In this case we say $l$ is of the ``second type''.
In both cases, we set $S_{l+1} = WS_l$.

Note that $j>s$ by virtue of the fact that
$w_r(D^{-1}S_lAS_l^{-\sigma} - I_n) < \epsilon$. As a consequence,
$w_r(D^{-1}W^\sigma D - I_n) < w_r(W - I_n)$, since if the leading term
on the right side is contained in $T_j t^j$, that on the left is 
the corresponding entry of $D^{-1} T_j D t^{pj}$, which has smaller
$w_r$-norm by the choice of $s$ and the bound $j>s$.
In particular, we have $w_r(D^{-1}S_{l+1}AS_{l+1}^{-\sigma} - I_n) < \epsilon$,
so the argument can be continued.

We now show that the $S_l$ converge in the Fr\'echet topology on
$\GL_{\an,r}$. Before doing so, we construct an upper bound on the size
of the coefficients of $A$ of positive index that will be stable under
conjugation by $S_l$. Namely, put $A = \sum_i A_i t^i$.
For $i$ in the value group of $L$, let $f(i) =
\max_{j\leq i}\{1, |D^{-1} A_j|\}$.
Let $i_0$ be the largest value such that $f(i_0) = 1$, and
extend $f$ to a smooth increasing
function $f(i)$ on $\RR$ which is identically 1 for
$i \leq i_0$. Now define
\[
g(i) = \exp\left(\int_0^i \max\{\frac{\Delta}{t}, \frac{f'(t)}{f(t)} \}
\right) \, dt.
\]
Then it is easily verified that $g(i) \geq f(i)$ for all $i$,
$(\log_p g(i))/i \to 0$ as $i \to +\infty$, and
$g(pi) \geq p^\Delta g(i)$ for $i \geq i_0$.
(The second condition follows from the inequality
$g(i) \leq if(i)$ obtained by bounding $\max\{x,y\}$
by $x+y$.)

By construction, $|D^{-1} A_i| \leq g(i)$ for all $i \geq 0$. Moreover, because
$g(pi) \geq p^\Delta g(i)$ for all $i$, it is easily verified by induction
on $l$ that for $D^{-1} S_l A S_l^{-\sigma} = \sum_i B_{l,i} t^i$, we have
$|B_{l,i}| \leq g(i)$ for all $i \geq 0$ and $l \geq 0$.

Now observe that $(r-r_l) j_l$ takes any given value only finitely
many times. To wit, if there were a smallest value that occurred infinitely
often, the corresponding values of $j_l$ would form a decreasing
sequence (since the $r$-minimal terms of $D^{-1}S_l A S_l^{-\sigma} - I_n$
would not change between these values of $l$, so $r_l$ would decrease on
these values). But $j_l$ lies in the (discrete) value group of $L$,
contradiction. Thus $(r-r_l)j_l$ takes any given value finitely many
times, and so goes to infinity with $l$.

Since $rj_l > (r-r_l)j_l$, $j_l$ also goes to infinity. On the other hand,
since $r_lj_l \leq g(j_l)$ and $g(i)/i \to 0$ as $i \to \infty$, $r_l \to 0$
as $l \to \infty$. By the same argument as in the previous paragraph,
$(r'-r_l)j_l \to \infty$ for $0 < r' < r$. Therefore the $S_l$ form a 
Cauchy sequence in the Fr\'echet topology, and their limit $S$ exists.

We finally prove that $DSAS^{-\sigma}-I_n$ has no coefficients of norm
greater than or equal to 1. We first prove there are no coefficients
of norm greater than 1. Put 
$D^{-1} S_l A S_l^{-\sigma} = \sum_i B_{l,i} t^i$
and $DSAS^{-\sigma} = \sum_i B_i t^i$. If
$|W_i| > 1$ for some $i$, let $j$ be the smallest index maximizing
$(\log_p |B_i|)/i$. Then for some $L$,
the coefficient of $t^j$ in $DS_lAS_l^{-\sigma} - I_n$ has norm
$|W_j|$ for all $l \geq L$. In fact, we can choose $L$ large enough so that
$(\log_p |B_{l,i}|)/i$ is maximized at $i=j$ for all $l \geq L$, by 
finding $i'$ such that $g(j)/j \leq (\log_p |B_i|)/i$ for all $j \geq i'$,
and choosing $L$ so that $|B_{l,i}| = |B_i|$ for $i \leq i'$.
For $l$ sufficiently large, $l$ cannot be of the first type, since $j_l <
i$ for all $l$ of the first type. So there exists $l \geq L$ of the
second type. By construction, $r_l = (\log_p |B_{l,j}|)/j$ and
$|B_{l+1,j}| < p^{r_lj}$, contradiction.

We conclude that $|B_i| \leq 1$ for all $i$. The proof that $|B_i| < 1$
for $i > 0$ is analogous but simpler: assuming the contrary, let
$j$ be the smallest index such that $|B_j| = 1$. Then for $l$ sufficiently
large, $|B_{l,i}| = |B_i|$ for $i \leq j$; but then any such $l$ would be
of the first type with $j_l = j$, and we would have $|B_{l+1,j}| < 1$,
contradiction.

In conclusion, we have that $D^{-1}SAS^{-\sigma}$ has entries in $\GLcon$
and $|D^{-1}SAS^{-\sigma}-I_n|<1$, as desired.
\end{proof}

The following proposition was inspired by an observation of
Kevin Buzzard and Frank Calegari, in their work on the Newton polygons
of Hecke operators acting on spaces of overconvergent modular forms.
\begin{prop} \label{prop:genspec}
Let $M$ be an $F$-crystal over $\Gimm$. Suppose there exists
a basis through which $F$ acts by the matrix $A$, where $D^{-1}A \equiv I_n
\pmod{\pi}$ for some diagonal matrix $D$ over $\calO$.
Then the slopes of the (generic) Newton
polygon of $M$ equal the valuations of the diagonal entries of $D$.
\end{prop}
\begin{proof}
Without loss of generality, assume $D$ has entries in $\calO_0$.
Let $B = D^{-1}AD$.
We produce a sequence of matrices $\{U_l\}_{l=1}^\infty$ such that $U_1 = I_n$,
$U_{l+1} \equiv U_l \pmod{\pi^{l}}$ and $U_l^{-1} BU_l^\sigma D^{-1} \equiv
I_n\pmod{\pi^l}$;
the $p$-adic limit $U$ of the $U_l$ will satisfy $BU^\sigma = UD$,
proving the proposition. The conditions for $l=1$ are satisfied by
the assumption that $D^{-1}A \equiv I_n \pmod{\pi}$.

Suppose $U_l$ has been defined. Let $U_l^{-1} BU_l^\sigma D^{-1} =
I_n + V$. For each $i,j$, let $W_{ij}$ be a solution of the
equation $W_{ij} - D_{ii} W_{ij}^\sigma D_{jj}^{-1} = V_{ij}$. Then
$W$ and $DW^\sigma D^{-1}$ are both congruent to 0 modulo $\pi^l$.
Put $U_{l+1} = (I_n+W)U_l$; then
\begin{align*}
U_{l+1}^{-1} B U_l^\sigma D^{-1} &=
(I_n+W)^{-1}U_l^{-1} B U_l^\sigma (I_n + W)^\sigma D^{-1} \\
&= (I_n+W)^{-1} U_l^{-1} B U_l^\sigma D^{-1} (I_n + D W^\sigma D^{-1}) \\
&= (I_n+W)^{-1}(I_n + V)(I_n + DW^\sigma D^{-1}) \\
&\equiv I_n-W+V+DW^\sigma D^{-1} = I_n \pmod{\pi^{l+1}}.  
\end{align*}
Thus the conditions for $U_{l+1}$ are satisfied, and the proposition
follows.
\end{proof}

Putting everything together, we obtain the following result.
\begin{prop} \label{prop:lowslope}
Let $M$ be an $F$-crystal over $\Gancon$ admitting a basis of eigenvectors
over $\Galgancon$. Suppose the lowest special slope of $M$ is $s$ occurring
with multiplicity $m$. Then there exists a $\Gcon$-lattice $N$ of dimension
$m$ which is $F$-stable with all slopes equal to $s$.
\end{prop}
\begin{proof}
Let $A$ be the matrix through which $F$ acts on some basis of $M$.
By Theorem~\ref{thm:constant}, there exists a matrix $D$ over a suitable
extension of $\calO_0$ and
a matrix $U$ over $\Galgancon$ such that $AU^\sigma = UD$. 
By Lemma~\ref{lem:anfact1},
for some finite extension $L$ of $K$,
there is an $F$-stable
$\GLcon$-lattice $M_1$ of full rank within $M \otimes_{\Gancon} \GLancon$
on which the action of $F$ is given by a matrix of the form $DV$, with
$V$ invertible over $\GLcon$. By applying $\sigma$ repeatedly, we may reduce
to the case where $L$ is separable.
By Proposition~\ref{prop:genspec},
viewed as an $F$-crystal over $\GLcon$, $M_1$ has the same special and
generic Newton polygons. By \cite[Corollary~6.6]{bib:me5},
$M_1$ has a subcrystal
$M_2$ of dimension $m$ with all slopes equal to $s$.

Otherwise put, the $\Galgcon$-module spanned by the eigenvectors of
$M \otimes_{\Gancon} \Galgancon$ of lowest slope descends to $\GLcon$.
Since this module admits an action of $\Gal(\overline{K}/K)$, $M_2$ admits an
action of $\Gal(L/K)$. We may apply Galois descent to the finite
extension $\GLcon/\Gcon$ of discrete valuation rings to conclude that
$M_2 = N \otimes_{\Gcon} \Gancon$ for some $\Gcon$-lattice $N$, which
is necessarily $F$-stable. The property of being \'etale is stable under
finite separable extensions, so $N$ is also \'etale, as desired.
\end{proof}

\section{Proof of the main result}

In this section, we assemble everything into the proof of the main theorem.
This stage of the proof is the first and only time where the connection
structure is used.

Let $M$ be an $(F,\nabla)$-crystal over $\Gcon$.
We may assume, by extending $\calO$, that the value group of $\calO$
contains all of the special and generic slopes of $M$.
Let $N
\subseteq M \otimes_{\Gcon} \Gancon$ be the lowest slope piece
of $M$ as given by Proposition~\ref{prop:lowslope}.
We show that $N \otimes_{\Gcon} \Gancon$ is stable under $\nabla$.
Choose a basis for $N$ and extend it to a basis of $M \otimes_{\Gcon}
\Gancon$. On this basis, $F$ acts by a block matrix of the form
$\begin{pmatrix} X & Y \\ 0 & Z \end{pmatrix}$ with $X$ an invertible
matrix over $\Gcon$. Let $\nabla$ acts by the block matrix
$\begin{pmatrix} P & Q \\ R & S \end{pmatrix}$. The relation
$\nabla F = p F \nabla$ implies the matrix identity
\[
\begin{pmatrix} P & Q \\ R & S \end{pmatrix}
\begin{pmatrix} X & Y \\ 0 & Z \end{pmatrix}
+
\begin{pmatrix} X & Y \\ 0 & Z \end{pmatrix}'
=
p
\begin{pmatrix} A & B \\ 0 & D \end{pmatrix}
\begin{pmatrix} P & Q \\ R & S \end{pmatrix}^\sigma,
\]
the lower left corner of which yields
$RX = pZR^\sigma$. We can write $X = U^{-1} U^\sigma$ with
$U$ over $\Galgcon$ by \cite[Lemma~2.2]{bib:me4},
and $Z = V^{-1} D V^\sigma$ with $V$ over $\Galgancon$ and $D$
a scalar matrix over $\calO_0$ whose entries have norms less than 1
(because of the construction of $N$ as the lowest slope piece of $M$).
Consequently, $VRU^{-1} = pD(VRU^{-1})^\sigma$,
which implies $VRU^{-1} = 0$ and so $R=0$, proving the claim.

Since $N \otimes_{\Gcon} \Gancon$ is stable under $\nabla$,
we may apply  an observation of Berger \cite[Lemme~V.14]{berger}
to conclude that $N$ itself is stable under $\nabla$.
Now Tsuzuki's theorem on the potential triviality
of unit-root $F$-crystals \cite{bib:tsu1} implies that $M_1$
is quasi-unipotent. By induction on the dimension of $M$, $M$ is
quasi-unipotent, as desired.

\section*{Acknowledgments}
The author was supported by a National Science Foundation Postdoctoral
Fellowship. Thanks to the organizers of the Dwork Trimester 
in Padova for their
hospitality, and to Laurent Berger and Johan de~Jong for helpful
discussions.

\end{document}